\documentclass{amsart}

\usepackage[style=alphabetic,sorting=anyt,sortcites=true,backend=biber,maxbibnames=5,maxcitenames=5]{biblatex}
\usepackage{hyperref}
\usepackage{amsmath,amssymb,amsthm}
\usepackage[a4paper, total={5.5in, 8.5in}]{geometry}
\usepackage{cleveref}
\usepackage{tikz-cd}
\usepackage{mathtools}
\usepackage{bm}
\usepackage{enumitem}

\newtheorem{theorem}{Theorem}[section]
\newtheorem{prop}[theorem]{Proposition}
\newtheorem{lemma}[theorem]{Lemma}

\theoremstyle{definition}

\newtheorem{example}[theorem]{Example}

\theoremstyle{remark}
\newtheorem{remark}[theorem]{Remark}

\newcommand{\nor}{\mathrel{\lhd}}

\newcommand{\Z}{\mathbb{Z}}
\newcommand{\N}{\mathbb{N}}
\newcommand{\eq}[1]{\equiv #1\bmod 3}

\newcommand{\I}{\mathcal{I}}

\newcommand{\ful}{\N\cup\{\infty\}}
\newcommand{\RR}[1]{\mathcal{R}(#1)}
\newcommand{\R}[1]{R(#1)}
\newcommand{\ke}[1]{\mathrm{ker}(#1)}
\newcommand{\ID}[1]{\mathrm{Id}_{#1}}
\newcommand{\lk}{\ldots}
\newcommand{\p}{\varphi}

\newcommand{\sen}{\ZZ\rtimes_{A'}\Z}
\newcommand{\f}{\psi}
\newcommand{\zn}{Z(H_n)}

\newcommand{\ti}{\tilde}
\newcommand{\pro}{\prod_{i=1}^n}
\newcommand{\su}{\sum_{i=0}^{\frac{n-1}{2}}}
\newcommand{\chari}{p_{_A}(x)}
\newcommand{\charin}{p_{_{\ii A}}(x)}
\newcommand{\g}{\gamma}
\newcommand{\gm}{\Gamma}
\newcommand{\e}{\varepsilon}
\newcommand{\deli}{\delta^i}
\newcommand{\deln}{\delta^{n-i}}
\newcommand{\ci}{\circ}
\newcommand{\ai}{a_i}
\newcommand{\an}{a_{n-i}}
\newcommand{\al}{\alpha}

\newcommand{\ab}{{[E,E]}}
\newcommand{\lm}{\lambda}
\newcommand{\lmi}{\lambda_i}
\newcommand{\C}{\mathbb{C}}

\newcommand{\GLL}{\mathrm{GL}_3\left(\Z\right)}
\newcommand{\GLLL}{\mathrm{GL}_4\left(\Z\right)}
\newcommand{\GL}{\mathrm{GL}_2\left(\Z\right)}
\newcommand{\SL}{\mathrm{SL}_2\left(\Z\right)}
\newcommand{\PSL}{\mathrm{PSL}_2\left(\Z\right)}
\newcommand{\G}[1]{\mathrm{GL}_{#1}\left(\Z\right)}
\newcommand{\GC}{\mathrm{GL}_2\left(\C\right)}
\newcommand{\Aut}[1]{\mathrm{Aut}(#1)}

\newcommand{\SE}[1]{\mathrm{Spec}_{\mathrm{R}}(#1)}
\newcommand{\rec}{Reidemeister classes }
\newcommand{\res}{Reidemeister spectrum }
\newcommand{\ren}{Reidemeister number }
\newcommand{\rens}{Reidemeister numbers }
\newcommand{\red}{Reidemeister }
\newcommand{\fingertof}{finitely generated torsion-free }

\newcommand{\op}{$R_{\infty}$~property}

\newcommand{\ii}[1]{{#1}^{-1}}
\newcommand{\im}[1]{\mathrm{Im}(#1)}
\newcommand{\dd}[1]{\lvert\mathrm{det}(#1)\rvert}
\newcommand{\di}[1]{\lvert\mathrm{det}(I-#1)\rvert}
\newcommand{\de}[1]{\mathrm{det}(#1)}

\newcommand{\gn}[1]{\langle #1\rangle}
\newcommand{\re}[2]{[#1]_{#2}}

\newcommand{\rez}[1]{[#1]_{\p}}
\newcommand{\rezz}[2]{\left[#1\right]_{#2}}
\newcommand{\ov}[1]{\overline{#1}}
\newcommand{\ba}[1]{\bar{#1}}

\newcommand{\quadr}{4\N\cup\{\infty\}}
\newcommand{\triples}{3\N\cup\{\infty\}}
\newcommand{\even}{2\N\cup\{\infty\}}
\newcommand{\con}[2]{\ii{#1}\,{#2}\,{#1}}

\newcommand{\cont}[2]{{#1}\,{#2}\,\ii{#1}}
\newcommand{\sem}{{\Z^2\rtimes_A\Z}}
\newcommand{\semm}{{\Z^3\rtimes_A\Z}}
\newcommand{\semn}{{\Z^n\rtimes_A\Z}}
\newcommand{\sehm}{{H_n\rtimes_\f\Z}}
\newcommand{\ZZ}{{\mathbb{Z}^2}}
\newcommand{\ZZZ}{{\mathbb{Z}^3}}
\newcommand{\QQ}{{\mathbb{Q}^2}}
\newcommand{\T}[1]{\mathrm{Tr}(#1)}
\newcommand{\W}{W_{-1}}

\newcommand{\real}{\mathbb{R}}
\newcommand{\rn}{\real^n}
\newcommand{\af}{\mathrm{Aff}(\rn)}
\newcommand{\afv}{\mathrm{Aff}(\real^4)}
\newcommand{\gr}{\mathrm{GL}_n(\real)}

\begin{document}
\title[Reidemeister spectra for solvmanifolds]{Reidemeister spectra for solvmanifolds in low dimensions} 

\author[K. Dekimpe]{Karel Dekimpe}
\author[S. Tertooy]{Sam Tertooy}
\author[I. Van den Bussche]{Iris Van den Bussche}
\address{KU Leuven Campus Kulak Kortrijk\\
E. Sabbelaan 53\\
8500 Kortrijk\\
Belgium}
\email{karel.dekimpe@kuleuven.be}
\email{sam.tertooy@kuleuven.be}
\email{iris.vandenbussche@kuleuven.be}

\thanks{Research supported  by long term structural funding -- Methusalem grant of the Flemish Government}

\maketitle
\date{\today} 

	\begin{center}
	This is an Accepted Manuscript of an article published by the Juliusz Schauder Center for Nonlinear Studies in Topological Methods in Nonlinear Analysis on Jun 2019, available online:  \href{http://dx.doi.org/10.12775/TMNA.2019.012}{http://dx.doi.org/10.12775/TMNA.2019.012}.
\end{center}

\begin{abstract}
	The Reidemeister number of an endomorphism of a group is the number of
	twisted conjugacy classes determined by that endomorphism. The
	collection of all Reidemeister numbers of all automorphisms of a group
	$G$ is called the Reidemeister spectrum of $G$. In this paper, we
	determine the Reidemeister spectra of all fundamental groups of
	solvmanifolds up to Hirsch length 4.
\end{abstract}

\section{Introduction}
Let $G$ be a group and $\p:G\to G$ an endomorphism. 
Consider the following equivalence relation on $G$:
\[x\sim_{\p}y\text{ if and only if }\exists z\in G: x=zy\p(z)^{-1}.\]
The equivalence classes under $\sim_{\p}$ are the \emph{Reidemeister classes of $\mathrm{\p}$} or the \emph{$\p$-twisted conjugacy classes}. 
We denote the set of these equivalence classes by $\RR\p$.  
The number of equivalence classes is called the \emph{Reidemeister number of $\p$} and is denoted by $\R\p$. 
If $\RR\p$ is infinite, we write $\R\p=\infty$.
Subsequently, the \emph{Reidemeister spectrum of $G$} is defined as
\[\SE G:=\{\R\p\mid \p\in\Aut G\}.\]
If $\SE G=\ful$, the Reidemeister spectrum of $G$ is said to be \emph{full}.
If $\SE G=\{\infty\}$, we say that $G$ has the \emph{\op}.

Reidemeister numbers of morphisms of groups correspond to Reidemeister number of self-maps on  topological spaces. These latter play a crucial role in the Nielsen theory of  fixed points and periodic points of these maps. We refer to \cite{jm06-1}, \cite{jian83-1} and \cite{kian89-1} for details about this. The study of groups having (or not having) the $R_\infty$ property was initiated by 
A.~Fel'shtyn and R.~Hill in the mid 90's and is now a very active research topic (see \cite{fh99-1,ft15-1,fn16-1,ddp09-1,dg14-1,gw03-1,gw09-1,gw09-2,roma11-1} for some papers in this area). In \cite{ft15-1} it was conjectured that finitely generated, residually finite groups without the $R_\infty$ property must be virtually solvable.
Having this in mind and also the topological meaning of a Reidemeister number (as being related to the study of fixed points), we will focus on a special class of  solvable groups, namely the fundamental groups of solvmanifolds. In fact, 
for the rest of this paper, we will study these  groups up to dimension $4$.
It is well known \cite{ov93-1} that a group $E$ is the fundamental group of a compact solvmanifold if and only if $E$ is an extension of the form
\begin{equation}\label{basis-extension}
    \begin{tikzcd}
    1 \arrow{r} & N \arrow{r}& E \arrow{r}{p} & \Z^k \arrow{r} & 1
\end{tikzcd}
\end{equation}
where $N$ is a finitely generated torsion-free nilpotent group and $k$ is a non-negative integer ($E=N$ when $k=0$). The dimension of the corresponding compact solvmanifold is the same as  the Hirsch length of the group $E$.

The aim of this paper is to determine the Reidemeister spectrum of all such  $E$ of Hirsch length at most $4$.
It turns out that many of these groups satisfy the \op. 
The groups that do not satisfy the $R_\infty$ property nor have full spectrum (abelian groups), have 
either \res $\{2, \infty\}$, $\{4,\infty\}$, $\{8,\infty\}$, $\even$, $4\N\cup\{\infty\}$, $6\N\cup\{\infty\}$ or $8\N\cup\{\infty\}$.

This paper is organised in five sections. 
We begin the paper by recalling 
some formulas to compute \red numbers. 
In Section~\ref{sectie3}, 
we first focus on the special case where the group $E$ is nilpotent.
In Section~\ref{sectie4}, we consider the groups $E$ of Hirsch length 3; in Section~\ref{sectie5}, we consider the groups $E$ of Hirsch length 4. 

\section{Preliminaries}\label{sectie2}
When determining the \res of the group $E$, an extension of $N$ by $\Z^k$, it will be convenient that the subgroup $N$ is characteristic. In this situation, we can easily determine \rens using the following \textit{addition formula}:
\begin{lemma}[{\cite[Lemma~2.1]{gw03-1}}]\label{formule}
Let $E$ be an extension of some group $G$ by $\Z^k$. Let $s$ be any (set-theoretic) section, so we have the exact sequence
with $p\circ s$ the identity on $\Z^k$:
\begin{center}
    \begin{tikzcd}
    1 \arrow{r} & G \arrow{r}& E \arrow{r}{p} & \Z^k \arrow{r}\arrow[bend left=33]{l}{s} & 1
\end{tikzcd}
\end{center}
Define the function $\al:\Z^k\to \Aut G$ by $\,\al(z)(g)=\cont{s(z)}{g}$.

Let $\p$ be an endomorphism of $E$ such that $\p(G)\subseteq G$. 
Write $\p|_G=\p'$ and let $\ba \p$ denote the induced endomorphism on the quotient $\Z^k$. 
Suppose $\{z_i\mid i\in\I\}$ is a complete set of representatives for the \rec of $\ba\p$. 
Then \[\R\p=\sum_{i\in\I} \R{\al(z_i)\circ\p'}.\]
\end{lemma}
In particular, $\R\p=\infty$ whenever $\R{\ba\p}=\infty$ or $\R{\p'}=\infty$. In general, for any automorphism $\p:G\to G$ inducing $\ba \p:G/H\to G/H$,
it always holds that if $\R{\ba\p}=\infty$, then $\R\p=\infty$ as well. Indeed, the map $\hat \pi:\RR\p\to\RR{\ba\p}: \rez x\mapsto\rezz{\ba x}{\ba\p}$ is well-defined and surjective. Hence, if $H$ is a characteristic subgroup of $G$ and $G/H$ has the \op, then $G$ has the $R_\infty$ property as well.

Since  the groups $E$ we are interested in in this paper are built by repeated extensions of the groups $\Z^n$, the following well-known formula is pivotal:

\begin{lemma}
Suppose $\p:\Z^n\to\Z^n $ is multiplication by
$M\in\mathrm{GL}_n(\Z).$ Then the \rec of $\p$ are the cosets of $\ker{(I-M)}$ in $\Z^n$, that is, $\RR\p=\Z^n/\ker{(I-M)}$. Moreover, $\R\p=\di M$ if this is non-zero, and $\R\p=\infty$ otherwise. 
\end{lemma}

In the sequel, we will also use $R(M)$ to denote the Reidemeister number $R(\varphi)$ of the automorphism $\varphi$ which is multiplication by $M$.

\medskip

If the subgroup $N$ in \eqref{basis-extension} is not characteristic, we can sometimes compute \rens using the \textit{averaging formula}. 
We briefly explain the setting of this formula. 
We refer to \cite{deki17-1} for a more general introduction and to \cite{char86-1,szcz12-1,wolf77-1} for the proofs of the results mentioned.

Consider the group of \textit{affine transformations} $\af=\rn\rtimes\gr$ on $\rn$ with multiplication $(a,A)\cdot(b,B)=(a+Ab, AB).$
An \textit{$n$-dimensional Bieberbach group} is a torsion-free cocompact discrete subgroup of $\rn\rtimes C$, where 
$C$ is a maximal compact subgroup of $\gr$.
Equivalently, an $n$-dimensional Bieberbach group $\gm$ is a torsion-free subgroup of $\af$ such that its subgroup of pure translations $\gm\cap\rn$ has finite index in $\gm$ and is a uniform lattice of $\rn$. A lattice of $\rn$ is a discrete and cocompact subgroup of $\rn$ and hence is isomorphic to $\Z^n$. From the geometric point of view, one mostly chooses $C$ to be the orthogonal group $O(n)$. Then $\Gamma$ is a subgroup of the group of Euclidean motions and the quotient manifold $\Gamma\backslash \rn$ inherits the flat metric structure from Euclidean space. From the algebraic point of view however, it is often easier to adopt another point of view. After an inner conjugation of $\af$, we may also assume that  $\gm\cap\rn$ is not only isomorphic to $\Z^n$, but  really coincides with $\Z^n$. From now on, we will always assume that this is the case.
The condition $\gm\cap\rn=\Z^n$ implies that
\begin{itemize}
    \item any element $(a,A)\in \gm$ has linear part $A$ in $\G n$, 
    \item any two elements $(a,A), (b,A)\in\gm$ are equal modulo $\Z^n$.
\end{itemize}
Hence $F:=\{A\in\G n\mid \exists a\in\rn: (a,A)\in\gm\}\cong\gm/\Z^n$. 
Note that $F$ is finite by definition. We call $F$ the \textit{holonomy group} of $\gm$.

Let $\p$ be an automorphism of $\gm$. 
The second Bieberbach Theorem says that $\p$ must be conjugation with some element  $\af$, that is, $\exists (m,M) \in \af$ such that 
$\p(\g)=(m,M)\g(m,M)^{-1}$ for all $\g\in\gm$. 
In particular, 
$\p((z,I))=(Mz,I)$ 
for all $z\in\Z^n$, hence $\p$ restricts to the automorphism $M\in\G n$ on $\Z^n$. In \cite{ll09-1} (see also \cite{hlp12-1}) we can find the following result:
\begin{lemma}\label{averagingformumla}
Let $\gm$ be a Bieberbach group with holonomy $F$. Let $\p:\gm\to\gm$ be an automorphism defined by $\p(\g)=(m,M)\g(m,M)^{-1}.$ Then $\R\p=\dfrac{1}{|F|}\displaystyle\sum_{A\in F} \R{AM}.$
\end{lemma}

\section{The nilpotent case}\label{sectie3}

We start with the special case where $E$ is itself nilpotent, corresponding to the subclass of nilmanifolds. In this case, already much is known, so we quickly review the \red spectra of the \fingertof nilpotent groups of Hirsch length at most $4$.  
We present the results in order of increasing nilpotency degree.

\subsection{\texorpdfstring{Nilpotency degree $1$.}{u}}
The abelian groups we have to consider
are the groups $\Z$, $\ZZ$, $\ZZZ$ and $\Z^4$. 
It is easy to see that $\SE \Z=\{2,\infty\}$ and $\SE {\Z^n}=\ful$ for $n\geq 2$. 
See for instance \cite{roma11-1}.
    
\subsection{\texorpdfstring{Nilpotency degree $2$.}{u}}
The finitely generated $2$-step nilpotent groups of Hirsch length at most 4
are of the form $H_n$ and $H_n\times\Z$, where $H_n:=\gn{x,y,z\mid [z,x]=1, [z,y]=1, [y,x]=z^n}$ for all $n\in\N.$
Roman'kov showed \cite[Section~3]{roma11-1} that $\SE{H_n}=\even$ for $n=1$, but his argument goes through for general $n$. 

We next compute the \res of $H_n\times\Z.$ 
Take a generator $u$ of $\Z$.
Let $\p$ be an automorphism of $H_n\times \Z$,
and let $\p'$ and $\ba \p$ denote the induced automorphism on the center $\gn{z,u}$ and on the quotient $\gn{\ba x, \ba y}$, respectively. 
Further, denote by $M\in\GL$ the matrix representing $\ba \p$ relative to the basis $\{\ba x, \ba y\}$. It is easy to check that $\p(z)^n=[\p(y),\p(x)]=z^{n\det(M)}$. Hence $\p(z)=z^{\det(M)}$. So $\p'$ is represented by the matrix
\[N:=\begin{pmatrix}
    \de{M} & r\\
    0 & \e
\end{pmatrix}\in\GL\]
for some $\e\in\{\pm 1\} $ and $r\in\Z$. 
If $\de M=1$ or $\e=1$, then $\R {\p'}=\infty$, and if $\de M=\e=-1$, then $\R{\p'}=4$. Moreover, the addition formula (Lemma~\ref{formule}) simplifies to $\R\p=\R{\p'}\R{\ba\p}$. Hence $\R\p\in\quadr$.

Conversely, take $m\in\N$ and consider the automorphism
\[\p_{m}:H_n\times\Z\to H_n\times\Z:x^ay^bz^cu^d\mapsto y^a(xy^m)^bz^{-c}u^{-d}.\]
Note that $\p_m$ induces the automorphisms
$N:=-I$
and $M:=\left(\begin{smallmatrix}
    0&1\\1&m
\end{smallmatrix}\right)$
on $\gn{z,u}$ and on $\gn{\ba x, \ba y}$, respectively. Hence $\R{\p_m}=\R{N}\R{M}=4m$. We conclude that $\SE{H_n\times\Z}=\quadr$ for all $n\in\N.$

\subsection{\texorpdfstring{Nilpotency degree $3$.}{u}}
The $3$-step nilpotent groups
necessarily have Hirsch length $4$. 
Gon\c calves and Wong treated in  \cite[Example~5.2]{gw09-2} an example of such a group and showed that this group has the \op. But in fact, using an analogous argument, one can show that all finitely generated torsion-free nilpotent groups of class 3 and Hirsch length $4$ have the \op.

\section{The 3-dimensional case}\label{sectie4}
In this section, we determine the \res of all fundamental groups of solvmanifolds of dimension at most $3$.
In the sequel, we denote by $\Z^l\rtimes_A \Z$, $A\in {\rm GL}_l(\Z)$, the semidirect product in which the generator of $ \Z$ is acting via $A$ on $\Z^l$.   More generally, we will use the notation $G\rtimes_\psi \Z$ for a semidirect product where the action is determined by an automorphism $\psi \in \Aut{G}$.

\medskip

Recall that the groups  $E$ we are interested in fit in  an exact sequence
\begin{center}
    \begin{tikzcd}
	1 \arrow{r} & N \arrow{r}& E \arrow{r} & \Z^k \arrow{r} & 1
\end{tikzcd}
\end{center}
with $N$ finitely generated torsion-free nilpotent.

\medskip

Then $E\cong\Z$ if $h(E)=1$ and either $E\cong\ZZ$ or $E\cong \Z\rtimes_{-1}\Z$ if $h(E)=2$. 

When $h(E)=3$, either $E\cong \Z^3$ or $E\cong H_n$ if $k=0$, and $E\cong\ZZ\rtimes\Z$ if $k=1.$
Moreover, if $h(E)=3$ and $k=2$, the following lemma says that $E$ can be viewed as a semidirect product  $\ZZ\rtimes\Z$ (so as an extension with $N\cong \Z^2$ and $k=1$) as well:

\begin{lemma}\label{extomdraaien}
Let $E$ be an extension of $\Z$ by $\Z^k$ of Hirsch length $k+1$.
Then $E\cong N\rtimes\Z$ for some \fingertof  nilpotent group $N$ of Hirsch length $k$.
\end{lemma}
\begin{proof}
By assumption, the group $E$ fits in the exact sequence
\begin{center}
	\begin{tikzcd}
		1 \arrow{r} & \Z \arrow{r}& E \arrow{r}{p} & \Z^k \arrow{r} & 1.
	\end{tikzcd}
\end{center}
Take generators $x_1,\lk, x_k$ of $\Z^k$.
The action of any element of $\Z^k$ on $\Z$ is either multiplication by $1$ or $-1$, so
we may assume that all $x_2,\lk,x_k$ act trivially on $\Z$, that is, $\Z$  is contained in the center of $\ii p(\gn{x_2,\lk,x_k})$. 
Moreover, the quotient $\ii p(\gn{x_2,\lk,x_k})/\Z\cong\gn{x_2,\ldots,x_k}$ is abelian, hence $\ii p(\gn{x_2,\lk,x_k})$ is a finitely generated, torsion-free, 2-step nilpotent group.
Since $E/\ii p(\gn{x_2,\lk,x_k})\cong \Z$, we conclude that $E\cong \ii p(\gn{x_2,\lk,x_k})\rtimes\Z$, as desired.
\end{proof}

Gon\c calves and Wong already showed that $\Z\rtimes_{-1}\Z$ has the $R_\infty$ property \cite[Theorem~2.2]{gw09-2}, so 
it remains to study \res of $\ZZ\rtimes\Z$.

We start by elaborating Lemma~\ref{formule}. 
\begin{lemma}\label{spec}
Suppose the subgroup $\Z^n$ is characteristic in $\semn$. 
Then \[\SE{\semn}=\{\infty\}\cup\{\R M+\R{AM}\mid M\in \G n, MA=\ii AM\}.\]
\end{lemma}
\begin{proof}
Take a generator $t$ of the quotient $\Z$. Let $\p$ be an automorphism of $\semn$ and suppose $\p|_{\Z^n}$ is multiplication with $M\in\G n$. 
 Note that $\R\p$ is infinite if $\p$ induces the identity on the quotient $\Z$. 
If not, $\p$ induces the automorphism $\ba\p=-\ID\Z$ on $\Z$. 
Using the representatives $1$ and $t$ for $\RR{\ba\p}$, the addition formula implies that $\R\p=\R M+\R{AM}$. 

Moreover, one easily verifies that there exists an automorphism $\p$ of $\semn$ inducing $M\in\G n$ on $\Z^{n}$ and $-\ID\Z$ on $\Z$ if and only if $MA=\ii AM$. 
This concludes the proof.
\end{proof}

The following lemma asserts that for most $A$, the subgroup $\Z^n$ is characteristic and  Lemma~\ref{spec} applies. It is an effortless generalization of  \cite[Lemma~2.1]{gw12-1}.

\begin{lemma}
\label{char}
If $A$ does not have $1$ as eigenvalue, the subgroup $\Z^n$ of $\semn$ is characteristic. 
\end{lemma}

We distinguish cases based on the eigenvalues of $A$.
\subsection{\texorpdfstring{The matrix $A$ has no eigenvalue $1$ or $-1$}{º}}\label{realorcomplex}

In \cite{gw03-1}, Gon\c calves and Wong investigated the $R_\infty$ property of this group. 
They found that $\sem$ does not have the $R_\infty$ property if and only if $A\in\SL$ and there exists $M$ in $\SL$ satisfying $\T M=0$ and $MA=\ii AM$. 
Writing 
\[A=
\begin{pmatrix}
    a & b\\
    c &  d\\
\end{pmatrix}
\text{ and }
M=
\begin{pmatrix}
    m & n\\
    p & -m\\
\end{pmatrix},\]
they showed that $MA=\ii AM$ is equivalent to $(a-d)m+bp+cn=0$. 
Hence $\sem$ does not have the $R_\infty$ property if and only if the system
\begin{equation}\label{system}
    \left\{
\begin{aligned}
    &-m^2-np=1\\
    &(a-d)m+bp+cn=0
\end{aligned}\right.
\end{equation} 
has a solution $(m,n,p)$ in $\Z^3$. 
Furthermore, they showed that both $\R M$ and $\R{AM}$ equal $1+\det(M)=2$, hence:
\begin{prop}[\cite{gw03-1}]\label{0or4}
Let $A\in\GL$ have eigenvalues different from $\pm 1$. Then $\sem$ has \res $\{\infty\}$ or $\{4,\infty\}$. Moreover, if $\det(A)=-1,$ then $\sem$ has the \op.
\end{prop}

Actually, if $A$ has complex eigenvalues, system (\ref{system}) has no solutions, so we moreover have:
\begin{prop}\label{complex}
Let $A\in\GL$ have non-real complex eigenvalues. Then $\sem$ has the \op.
\end{prop}
\begin{proof}
Note that the eigenvalues of $A$ are the roots of the polynomial 
\begin{align*}
    \de{xI-A}
    &=x^2-\T A x+1.
\end{align*} 
When $A$ has complex eigenvalues, we therefore must have $\T A^2-4<0.$ 

Suppose for a contradiction that $\sem$ does not have the \op, that is, the system (\ref{system}) has an integral solution $(m,n,p)$. 
Note that $c\neq 0$, for otherwise $A$ would have eigenvalues $\pm1$. 
Substituting the second equation in the top one gives 
\[-m^2+\dfrac{a-d}{c}mp+\dfrac{b}{c}p^2=1.\]
Note that $p\neq 0$. Denoting $y:=m/p\in\mathbb{Q}$, this relation rewrites to 
\[-y^2+\dfrac{a-d}{c}y+\dfrac{b}{c}-\dfrac{1}{p^2}=0.\]
In particular, the discriminant $\bigl((a-d)/c\bigr)^2+4(b/c-1/p^2)$ must be positive. 
Using $ad-bc=\de A=1$, this simplifies to
\begin{align*}
    0\leq\left(\dfrac{a-d}{c}\right)^2+4\left(\dfrac{b}{c}-\dfrac{1}{p^2}\right)
    &=\dfrac{a^2+d^2-2ad}{c^2}+\dfrac{4b}{c}-\dfrac{4}{p^2}\\
    &=\dfrac{(a^2+d^2-2ad+4bc)p^2-4c^2}{c^2p^2}\\
    &=\dfrac{(\T A^2-4)p^2-4c^2}{c^2p^2}\cdot
\end{align*}
Since $\T A^2-4<0$, we have reached a contradiction, and the proof is complete.
\end{proof}

In general, determining which $A$ allow an integral solution to (\ref{system}) is hard and we will not pursue this further. Instead, we proceed with the remaining case where $A$ has $1$ or $-1$ as eigenvalue. 
Note that this implies that both eigenvalues are $\pm 1$.
\subsection{\texorpdfstring{The matrix $A$ has eigenvalues $1$ or $-1$.}{u}}

Let $\e$, $\delta\in\{\pm 1\}$ be the eigenvalues of $A$. 
Choose an eigenvector $v\in\QQ$ corresponding to $\e$. 
Clearing denominators if necessary, we may assume that $v\in\ZZ$ and that $v$ extends  to a basis $\{v,w\}$ of $\ZZ$. So we can write 
\[A=\begin{pmatrix}
    \e & r\\
    0 &  \delta\\
\end{pmatrix}
\]
for some integer $r$. 

Let $W_\e:=\{z\in\ZZ\mid A(z)=\e z\}$ denote the eigenspace of $\e$. 
We make the following easy observations, the proof of which we leave to the reader:
\begin{lemma}\label{chareig}
Let $A\in\G n.$
\begin{enumerate}
    \item If $A$ does not have $1$ as eigenvalue, the eigenspace $W_{-1}$ is characteristic in $\semn.$
    \item The group $\semn$ has center $W_1\times \gn {t^d}$, where $d$ is the order of $A$ if this order is finite. If $A$ has infinite order, we set $d:=0$.
\end{enumerate}
\end{lemma}

Swapping $\e$ and $\delta$ if necessary, we have the following three cases.

\medskip\noindent\textit{Case 1: $\e=\delta=1$.} 
One easily sees that in this case, $\sem$ is nilpotent. Specifically, $\sem\cong\ZZZ$ if $A=I$ and $\sem\cong H_n$ for some $n\in\N$ if $A\neq I$. 

\medskip\noindent\textit{Case 2: $\e=\delta=-1$.}
We further distinguish the cases $r\neq0$ and $r= 0$.
\begin{prop}\label{2x-1}
Let $A\neq -I$ in $\GL$ have repeated eigenvalue $-1$. Then $\sem$ has the \op.
\end{prop}
\begin{proof}
The eigenspace of $-1$ equals $W_{-1}=\langle v\rangle$. 
Moreover, the quotient $(\sem)/W_{-1}\cong \Z\rtimes_{-1}\Z$ has the $R_\infty$ property \cite[Theorem~2.2]{gw09-2}. Hence 
$\sem$ has the $R_\infty$ property as well.
\end{proof}

\begin{prop}\label{-i2}
The group $\Z^n\rtimes_{-I}\Z$ has Reidemeister spectrum $2\N\cup\{\infty\}$ for all $n\geq2.$
\end{prop}
\begin{proof}
Note that any $M\in\G n$ satisfies $(-I)M(-I)=M$. Consequently, Lemma~\ref{spec} implies that
\[\SE{\Z^n\rtimes_{-I}\Z}=\{\infty\}\cup\{\R M+\R{-M}\mid M\in \G n\}.\]
Recall that $\R M+\R{-M}$, if finite, equals $\di M+\dd{I+M}$. 
Since $-1\equiv1\bmod 2$, surely $\R M+\R{-M}\in\even$.

Conversely, for $m$ in $\N\cup\{0\}$, consider the $(n\times n)$-matrix
\begin{align*}
 M_m:=
\left(
\begin{array}{c|c}
  \begin{array}{ccc} 0   & \cdots & 0 \end{array} & 1 \\
  \hline
  I_{n-1}
  & \begin{array}{c}0\\ \vdots\\ 0\\ m\end{array}
\end{array}
\right).   
\end{align*}
One easily verifies that $\di{M_m}+\dd{I+M_m}$  equals $2m$ if $n$ is even and $2(m+1) $ if $n$ is odd. This completes the proof.
\end{proof}

\medskip\noindent\textit{Case 3: $\e=1$, $\delta=-1$.}

\begin{prop}\label{1en-1}
Let $A$ in $\GL$ have eigenvalues $1$ and  $-1$. Then $\sem$ has the \op.
\end{prop}
\begin{proof}
The center of $\sem$ equals $Z=\langle v\rangle\times \langle t^2\rangle$. 
The quotient $(\sem)/Z\cong\Z\rtimes_{-1}\Z_2$ satisfies the $R_\infty$ property \cite[Proposition~2.3]{gw09-2}, hence so does $\sem$.
\end{proof}

\subsection{Conclusion}

We summarise our findings in 
the following table:

\[\begin{array}{cl|c}
\hline
\multicolumn{2}{c|}{A}&\SE \sem \\[1em]
\text{eigenvalues }&\\ \hline\hline
-1,-1&
\begin{array}{l}A=-I\\A\neq -I\end{array}
&
\begin{array}{c}\even\\ \{\infty\}\end{array}
\\[1em] 
1,-1& &\{\infty\}
\\[0.5em]
\text{real}\neq\pm1
&
\begin{array}{l}\det(A)=1\\\det(A)=-1\end{array}
&
\begin{array}{c}\{\infty\}\text{ or }\{4,\infty\}\\ \{\infty\}\end{array}
\\[1em] 
\text{non-real}& &\{\infty\}\\\hline
\end{array}\]

\section{The 4-dimensional case}\label{sectie5}

In this section, we consider extensions $E$ as in \eqref{basis-extension} where $h(E)=4$.

Again, the special case $k=0$ corresponds to the nilpotent groups of Section~\ref{sectie3}, and 
by Lemma~\ref{extomdraaien}, all other groups will appear in the situations $k=1$ and $k=2$. 

In case $k=1$, the sequence splits, hence $E\cong \ZZZ\rtimes \Z$ or $E\cong H_n\rtimes\Z$.

In case $k=2$, the group $E$ is an extension of $\ZZ$ by $\ZZ$:
\begin{center}
	\begin{tikzcd}
		1 \arrow{r} & \ZZ \arrow{r}& E \arrow{r}{p} & \ZZ\arrow{r}  & 1.
	\end{tikzcd}
\end{center}
Let $\al:\ZZ\to \GL$ denote the induced action of the quotient $\ZZ$ on $\ZZ$. 
Take  generators $x$, $y$ of this quotient, write $\al(x)=A$ and $\al(y)=B$. 
As $x$ and $y$ commute, the matrices $A$ and $B$ commute as well. 
This will severely limit the possible values of $A$ and $B$. 
We need the following lemma:
\begin{lemma}\label{com}
Let $M\neq \pm I$ in $\GL$ have finite order. 
Then the centraliser of $M$ in $\GL$ is the subgroup $
\langle -I, M \rangle$.
\end{lemma}
\begin{proof}
This lemma is easily checked by a case-by-case study, using the fact that there are, up to conjugacy, only six finite cyclic subgroups in $\GL$ (see e.g.\ \cite[page 179]{newm72-1}).
\end{proof}

Using this, we will show that there are essentially three possibilities for $A$ and $B$:

\begin{lemma}\label{gener}
We can choose $x$ and $y$ such that one of the following conditions holds:
\begin{enumerate}
    \item $B=I$,
    \item $B=-I$ and $A$ has infinite order,
    \item $B=-I$ and $A\neq \pm I$ has order $2$.
\end{enumerate}
\end{lemma}
\begin{proof}
Note that $\al:\ZZ\to\GL$ is not injective, as $\PSL\cong\Z_2\ast\Z_3$ does not contain $\ZZ$ so neither does $\GL$.
Thus either $\ZZ/\ke \al$ is finite (\textbf{Case 1}) or $\ZZ/\ke \al\cong\Z\oplus\Z_p$ for some $p\in\N$ (\textbf{Case 2}). 

\medskip\noindent\textbf{Case 1.} If $\ZZ/\ke \al$ is finite, both $A$ and $B$ have finite order. 
If $A=I$ or $B=I$, we are in situation~(1). Moreover, if both $A=-I$ and $B=-I$, replacing $y$ by $xy$ results in situation~(1) as well. Hence, it remains to prove the lemma when $A$ or $B$ is different from $\pm I$. 
We may assume $A\neq \pm I$. By Lemma~\ref{com}, there exists $k\in\Z$ such that $B=\pm A^k$. 
Replacing $y$ by $x^{-k}y$ if necessary, we can further reduce to $B=\pm I$. 
In case $B=I$, we are back in situation~(1), so assume $B=-I$. 
If $A$ happens to have order 2, we are in situation~(3). 
Otherwise, $A$ has order $3$, $4$ or $6$.
If $A$ has order $3$, the order of $-A$ is $6$. Replacing $x$ by $xy$ if necessary, we are left with the situation where $A$ has order $d\in\{4,6\}$. 
However, then $A^{d/2}=-I$ and replacing $y$ by $x^{d/2}y$ delivers situation~(1).

\medskip\noindent\textbf{Case 2.} Suppose $\ZZ/\ke \al\cong\Z\oplus\Z_p$. 
Take $u\in \ZZ$ such that $\ba u$ generates $\Z_p$. 
If we choose $x$, $y$ such that $y^k=u$ for some $k$ in $\Z$, 
then $B$ has finite order, whereas $A$ must have infinite order as $\ZZ/\ke\al$ is infinite. 
By Lemma~\ref{com}, this only happens when $B=\pm I$, corresponding to situations~(1) and (2). 
This completes the proof.
\end{proof}

So, we may assume that $x$, $y$ are generators as in Lemma~\ref{gener} above. 
Choose a section $s:\ZZ\to E$ and write $u=s(x)$ and $t=s(y)$. Then
$E\cong (\ZZ\rtimes_B\Z)\rtimes_\f\Z$ with action
\[\f(zt^l)=uzt^lu^{-1}=uzu^{-1}(utu^{-1})^l=A(z)(n_0t)^l,\quad z\in\ZZ,\,l\in\Z,\]
setting $n_0:=utu^{-1}t^{-1}\in\ZZ$. Thus either $E$ is isomorphic to $\Z^3\rtimes \Z$, 
or $E$ is isomorphic to $(\ZZ\rtimes_{-I}\Z)\rtimes_\f\Z$ where $\f|_\ZZ=A$ has infinite order or $A\neq\pm I$ and $A^2=I$.

Hence, it remains to study the Reidemeister spectrum of the groups $\Z^3\rtimes \Z$, $H_n\rtimes\Z$ and $(\ZZ\rtimes_{-I}\Z)\rtimes_\f\Z$.

\subsection{\texorpdfstring{The semidirect product $\Z^3\rtimes_A \Z$}{y}}\label{semm}

In this subsection, we determine the \res of the group $\semm$. 
We distinguish cases based on the eigenvalues of $A$.

\subsubsection{\texorpdfstring{The matrix $A$ does not have eigenvalue $1$}{u}}\label{semi3}

According to Lemma~\ref{spec}, we must determine which $M$ satisfy $MA=\ii A M$. 
If any such $M$ exists, $A$ is conjugate to its inverse.
The following lemma says that this only happens if $A$ has eigenvalue~$-1$, as we assume that 1 is not an eigenvalue of $A$.

\begin{lemma}\label{eig-1dekimpe}
Let $n\in\N$ be odd, and $A\in\mathrm{GL}_n(\Z)$. 
Suppose $A$ is conjugate to $\ii A$. 
Then $A$ has eigenvalue $\de A$.
\end{lemma}
\begin{proof}
Set $\delta=\de A$ and write $\lm_1,\lk,\lm_n$ the eigenvalues of $A$. 
Let $p_{_A}$ and $p_{_{\ii A}}$ denote the characteristic polynomial of $A$ and $\ii A$, respectively. 
Since $A$ and $\ii A$ are conjugate, their characteristic polynomials coincide, hence
\begin{align*}
    \chari=\charin&=\pro (x-\ii\lmi)\\
    &=\dfrac{x^n}{\lm_1\lk\lm_n}\pro (\lmi-\ii x)\\
    &=(-1)^n\delta x^n\pro (\ii x-\lmi)\\
    &=-\delta x^np_{_A}(\ii x).
\end{align*}
Writing $\chari=\sum_{i=0}^n \ai x^i$, this implies that $\ai=-\delta\an$, for $x^np_{_A}(\ii x)$ equals the reciprocal polynomial $\sum_{i=0}^n\an x^i$ of $p_{_A}$.
It easily follows that
\[p_{_A}(\delta)=\sum_{i=0}^n \ai\deli=\su (\ai\deli+\an\deln)=\su \deli(\ai+\delta\an)=0,\]
that is, $A$ has eigenvalue $\delta$.
\end{proof}

Subsequently, if neither $1$ nor $-1$ are eigenvalues of $A$, 
Lemma~\ref{spec} immediately implies the following:

\begin{prop}\label{differentfromApm1}
Let $A\in\GLL$ have eigenvalues different from $1$ and $-1$. Then $\semm$ has the \op.
\end{prop}

So it remains to study the $A$ that do have eigenvalue $-1$. 
Recall from Proposition~\ref{-i2} that $\SE\semm=\even$ if $A=-I$.
If $A\neq -I$, we have:

\begin{prop}\label{completespace-1}
Let $A\in\GLL$ be different from $-I$ and suppose that $A$ has eigenvalue~$-1$, but not eigenvalue~$1$. 
Then $\semm$ has the \op.
\end{prop}
\begin{proof}
Let $\{u,v,w\}$ be a basis of $\ZZZ$ transforming $A$ into  a matrix of the form \[\begin{pmatrix}
    -1 & C\\ 0 & A'
\end{pmatrix}\]
for some $A'\in\GL$, $C\in\Z^{1\times 2}.$ 
 Suppose first that $-1$ is not an eigenvalue of $A'$.
Let $\p$ be an automorphism inducing
\[
M:=
\begin{pmatrix}
    m & N\\ P & Q
\end{pmatrix}\in\GLL
\]
on $\ZZZ$
for some $m\in \Z$, $N\in\Z^{1\times 2}$, $P\in\Z^{2\times 1}$ and $Q\in\Z^{2\times 2}$.
Since $W_{-1}$ is characteristic, we must have
$P=0$, implying $m\in\{\pm1\}$.
As a result, either $M$ or $AM$ has eigenvalue $1$, or equivalently, either $\R M$ or $\R{AM}$ is infinite. 
Hence $\semm$ has the \op.

 If $A'$ does have eigenvalue $-1$, then $-1$ is the only eigenvalue of $A'.$
If $W_{-1}=\gn u$, it is easy to check that $A'\neq -I$. The quotient $\semm/\W\cong \ZZ\rtimes_{A'}\Z$ has the $R_\infty$ property by Lemma~\ref{2x-1}, hence so has $\semm$. If $W_{-1}\neq\gn u$, then $\semm/\W\cong\Z\rtimes_{-1}\Z$ has the $R_\infty$~property as well. This completes the proof.
\end{proof}

\subsubsection{The matrix $A$ has eigenvalue $1$}

We present our results in decreasing order of the 
algebraic multiplicity of  the eigenvalue $1$.

If $1$ is the only eigenvalue of $A$, the group $\semm$ is nilpotent.

If the eigenvalue $1$ has algebraic multiplicity two, we have the following:

\begin{prop}\label{alg2}
Suppose $1$ is an eigenvalue of $A\in\GLL$ of algebraic multiplicity two.
Then $\semm$ has the \op.
\end{prop}
\begin{proof}
If $A$ has eigenvalue $1$ of algebraic multiplicity two, the remaining eigenvalue must be $-1$. 
Accordingly, let $\{u,v,w\}$ be a basis of $\ZZZ$ transforming $A$ into  a matrix of the form
\[\begin{pmatrix}
    1 & r & s\\ 0& 1&n\\0&0&-1
\end{pmatrix}\]
for some $r,s,n\in\Z$. If $r$ is zero, $A$ has order two, so the center of $\semm$ is generated by $u$, $v$ and $t^2$. Moreover, the quotient $\semm/\gn{u,v,t^2}\cong\Z\rtimes_{-1}\Z_2$ has the \op, hence so has $\semm$. 

Similarly, if $r$ is nonzero, the order of $A$ is infinite. 
Thus $\semm$ has center $\gn u$ and $\semm/\gn u$ has the $R_\infty$~property by Lemma~\ref{1en-1}. 
We conclude that $\semm$ has the $R_\infty$~property as well, as desired.
\end{proof}

If the eigenvalue $1$ is not repeated, let $\{x,y,z\}$ be a basis of $\ZZZ$ such that $A$ takes the form
\begin{equation}\label{a'}
    \begin{pmatrix}
    1 & \begin{matrix}r&s\,\end{matrix}\\
    \begin{matrix}0\\0\end{matrix} & A'
\end{pmatrix}
\end{equation}
for some $r,s\in\Z$ and $A'\in\GL$ not having eigenvalue $1$. 
Thus $A'$ has either repeated eigenvalue $-1$ or real or complex eigenvalues. 
This leaves the following possibilities:
\begin{enumerate}[label=(\roman*)]
    \item $A'=-I$,
    \item $A'=\left(\begin{smallmatrix}-1&n\\0&-1
    \end{smallmatrix}\right)$ for some $n\neq 0$,
    \item $A'$ has real eigenvalues not equal to $\pm 1$,
    \item $A'$ has non-real complex eigenvalues.
\end{enumerate}
In cases (ii) and (iii), the matrix $A$ has infinite order, hence $x$ generates the center of $E=\semm$. 
Moreover, in $E/\gn x\cong \ZZ\rtimes_{A'}\Z$, the subgroup $\ZZ=\gn{\ba y,\ba z}$ is characteristic by Lemma~\ref{char}, hence $\ZZZ$ is characteristic in $E$.
Consequently, we can apply Lemma~\ref{spec} to compute the \res of $E$. 
To this end, take $M\in\GLL$ satisfying $MA=\ii AM$.
Since $\gn x=Z(E)$ is characteristic, we can write
\[M=\begin{pmatrix}
    m & N\\ 0 & Q
\end{pmatrix}\]
for some $m\in \{\pm1\}$, $N\in\Z^{1\times 2}$ and $Q\in\Z^{2\times 2}$. If $m=1$, clearly $\R M=\infty$. 
So assume $m=-1$. 
Imposing $AMA=M$ explicitly gives
\begin{equation*}    
\begin{pmatrix}
    -1 & -C+NA'+CQA'\\0 &A'QA'
\end{pmatrix}=
\begin{pmatrix}
        -1 & N\\ 0 & Q
\end{pmatrix},
\end{equation*}
hence $A'QA'=Q$. 
Moreover, it is easy to see that $\R M+\R{AM}=2(\R Q+\R{A'Q}).$ 
Thus $\SE\semm\subseteq 2\SE\sen$. 
Therefore, Lemma~\ref{2x-1} and \ref{0or4} immediately imply the following results.
\begin{prop}\label{nonzeron}
Suppose $A$ is conjugate to a matrix of the form (\ref{a'}) where $A'=\left(\begin{smallmatrix}-1&n\\0&-1
    \end{smallmatrix}\right)$ for some $n\neq 0$. Then $\semm$ has the \op.
\end{prop}
\begin{prop}\label{0or8}
Suppose $A$ is conjugate to a matrix of the form (\ref{a'}), where $A'$ has real eigenvalues $\neq \pm 1$.
Then $\semm$ has \res $\{\infty\}$ or $\{8,\infty\}$. 
Moreover, if $\ZZ\rtimes_{A'}\Z$ has the \op, so has $\semm$. 
\end{prop}

\begin{remark}
In general, the converse statement to Proposition~\ref{0or8} is not true: there exists $A$  
such that $\semm$ has the \op, but $\sen$ does not. 
Indeed, consider
\[A=\begin{pmatrix}
        1 & 0 & 1\\ 0 &5 & 2\\ 0&2&1
\end{pmatrix},\text{ so that }A':=\begin{pmatrix}
        5&2\\2&1
\end{pmatrix}.\]
In the above argument, $\SE\semm=\{8,\infty\}$ precisely when there exists $Q$ in $\SL$ satisfying $A'QA'=Q$ and 
\begin{equation}\label{deconditie}
    -N=C(I-QA')\ii{(I-A')}\in\Z^{1\times 2},
\end{equation} where $C=(0,1)$.
A straightforward calculation shows that if $Q\in \SL$ satisfies $A'QA'=Q$, all entries of $I-QA'$ are odd. 
At the same time, it is easy to check that 
\[\ii{(I-A')}=
\dfrac{1}{2}\begin{pmatrix}
        0&-1\\-1&2
\end{pmatrix}.\]
Subsequently, condition~(\ref{deconditie}) never holds. 
At the same time, $Q' = \left(\begin{smallmatrix}
        0 & 1\\ -1  &0 \end{smallmatrix}\right)$ satisfies $AQ'A=Q'$ and $R(Q')+R(AQ')=4$. By Lemma~\ref{spec}, 
we conclude that the group $\semm$ has the \op, while $\sen$ does not.
\end{remark}

In cases~(i) and (iv), the matrix $A$ has finite order, say order $d$. 
In this situation, we can use the averaging formula, Lemma~\ref{averagingformumla}, to compute \rens on $\semm$. 
Indeed, setting 
$\ti A:=
\left(\begin{smallmatrix}
1&0\\0&A
\end{smallmatrix}\right)\in\GLLL,$
the embedding 
\[i:\semm\to \afv:
x^ay^bz^ct^n\mapsto 
\left(
\begin{pmatrix}n/d \\ a\\b\\c\end{pmatrix} ,\,\ti A^n
\right)
\]
identifies $\semm$ with 
a $4$-dimensional Bieberbach group with holonomy  $\{I,\ti A,\lk,\ti A^{d-1}\}$.
Let $\p$  be an automorphism of $\semm$. 
Let $N$ denote the matrix representing $\p$ on the center $\gn{t^d,x}$ 
and let $M$ denote the matrix representing the induced automorphism $\ba\p$ on  $\ZZ=\gn{\ba y ,\ba z}$. 
The restriction of $\p$ to $\gn{t^d,x,y,z}\cong\Z^4$ is thus represented by a matrix of the form
\[\ti M:=
\begin{pmatrix}
N&*\\0&M
\end{pmatrix}\in\GLLL.\]
Then for all integers $i$, the matrix $\ti A^i\ti M$ has the form 
\[\ti A^i\ti M
=
\begin{pmatrix}
I&*\\0&{A'}^i
\end{pmatrix}
\begin{pmatrix}
N&*\\0&M
\end{pmatrix}
=
\begin{pmatrix}
N&*\\0&{A'}^iM
\end{pmatrix},
\]
showing $\R{\ti A^i\ti M}=\R N\R{ {A'}^i M}$. The averaging formula hence says that 
\begin{align*}
    \R\p
    &=\dfrac{1}{d}\R N\sum_{i=0}^{d-1} \R{ {A'}^i M}.
\end{align*}

We start with case (i), so that $A$ has order 2. Tahara showed \cite[Proposition~2]{t71-1} that in this case, $A$ is conjugate over $\GLL$ to the matrix
\begin{equation}\label{0of1}
    \begin{pmatrix}
    1 & 0 & \delta\\ 0& -1&0\\0&0&-1
\end{pmatrix},
\end{equation}
where $\delta$ is either $0$ or $1$. 
\begin{prop}\label{delta0}
Suppose $A$ is of the form (\ref{0of1}) where $\delta\in\{0,1\}$. 
\begin{enumerate}
    \item If $\delta=0$, the group $\semm$ has \res $2\N\cup\{\infty\}.$
    \item If $\delta=1$, the group $\semm$ has \res $4\N\cup\{\infty\}.$
\end{enumerate}
\end{prop}
\begin{proof}
Write $E=\semm$. 
Take an automorphism $\p$ of $E$. 
The center of $E$ is generated by $x$ and $t^2$, so we can write
\[\p(x)=x^at^{2b},\; \p(y)=x^cy^dz^et^f, \;\p(z)=x^gy^hz^it^j,\; \p(t)=x^ky^lz^mt^n\]
for some $a,\lk,n$ in $\Z$. 
Moreover, $f$ and $j$ are even and $n$ is odd as $\ZZZ\times \gn{t^2}$ is characteristic in $E$.

The map $\p$ is a morphism if and only if it respects the relations
\[[x,y]=1,\; [x,z]=1,\; [y,z]=1,\; [x,t]=1,\;[y,t]=y^{-2},\;[z,t]=x^\delta z^{-2}.\]
As both $f$ and $j$ are even, $\p$ already respects the first four relations, so we only need to examine the last two relations. Equating $\p(y)^{-2}=[\p(y),\p(t)]$ gives
\begin{align*}
x^{-2c}y^{-2d}z^{-2e}t^{-2f}&=t^{-f}z^{-e}y^{-d}x^{-c}\,t^{-n}z^{-m}y^{-l}x^{-k}\,x^cy^dz^et^f\,x^ky^lz^mt^n\\
&=z^{-e}y^{-d}y^{-d}x^{\delta e}z^{-e}t^{-n}t^n\\
&=x^{\delta e}y^{-2d}z^{-2e},\\
\intertext{so $\delta e=-2c$ and $f=0$. Similarly, equating $\p(x)^\delta\p(z)^{-2}=[\p(z),\p(t)]$ gives}
x^{\delta a-2g}y^{-2h}z^{-2i}t^{2\delta b-2j}&=t^{-j}z^{-i}y^{-h}x^{-g}\,t^{-n}z^{-m}y^{-l}x^{-k}\,x^gy^hz^it^j\,x^ky^lz^mt^n\\
&=z^{-i}y^{-h}y^{-h}x^{\delta i}z^{-i}t^{-n}t^n\\
&=x^{\delta i}y^{-2h}z^{-2i},
\end{align*} 
hence $\delta i=\delta a-2g$ and $\delta b=j$.

The averaging formula says that $\R\p=\dfrac12 \R N(\R M+\R{-M})$ where
\[N:=\begin{pmatrix}
    n&b\\2k+\delta m&a
\end{pmatrix}\text{ and }
M:=\begin{pmatrix}
     d&h\\e&i
 \end{pmatrix}.
\]
We distinguish two cases.

\medskip\noindent\textbf{Case 1.} Suppose first that $\delta=1$.
Then $e=-2c$ is even, $a\equiv i \bmod 2$ and $b=j$ is even. 
As $\R M=\pm 1$ is odd and $e$ is even, both $d$ and $i$ (and thus also $a$) are odd.
So $b$ and $e$ are even and $n$, $a$, $d$ and $i$ are odd.
It follows that  $\R N$ is either infinite or even. In addition, if $\R M + \R{-M}$ is finite, then
\[\R M+\R{-M}=|1+\det(M)-\T M|+|1+\det(M)+\T {M}|\in 4\N,\] as $\T M=d+i$ is even. 
Hence $\SE E\subseteq 4\N\cup\{\infty\}$.

Conversely, for $\al\in\N$, setting 
\[
\p_\al(x)=x^{1-2\al}t^{4\al},\; 
\p_\al(y)=x^{-1}y^{-1}z^2, \;
\p_\al(z)=x^{1-\al}yz^{-1}t^{2\al},\; 
\p_\al(t)=zt^{-1}
\] defines an automorphism $\p_\al$ of $E$
inducing the automorphisms
\[
\begin{pmatrix}
-1&2\al\\1&1-2\al 
\end{pmatrix}
\text{ and }
\begin{pmatrix}
    -1&1\\2&-1
\end{pmatrix}
\in\GL\]
on $Z(E)=\gn{ t^2,x}$ and $\ZZ=\gn{\ba y, \ba z}$, respectively. 
It is easy to check that $\R{\p_\al}=
4\al,$ hence in this case, $\semm$ has \res $\quadr$.

\medskip\noindent\textbf{Case 2.} If $\delta=0$,
then $\R N$ is either infinite or even as $n$ is odd. Moreover, $\R M+\R{-M}$, if finite, is always even, hence $\R\p\in\even $.

Conversely, for $\al\in\N$, setting 
\[
\p_\al(x)=x^{1-2\al}t^{2},\; 
\p_\al(y)=z, \;
\p_\al(z)=yz,\; 
\p_\al(t)=x^{\al}t^{-1}
\] defines an automorphism $\p_\al$ of $E$ inducing 
\[\begin{pmatrix}
-1&1\\2\al&1-2\al 
\end{pmatrix}
\text{ and }
\begin{pmatrix}
    0&1\\1&1
\end{pmatrix}
\in\GL\]
on $Z(E)=\gn{ t^2,x}$ and $\ZZ=\gn{\ba y, \ba z}$, respectively. 
It is easy to check that $\R{\p_\al}=2\al,$ which concludes the proof. 
\end{proof}

If $A'$ has complex eigenvalues, 
the order $d$ of $A'$ is either $3$, $4$ or~$6$. Additionally, if $d\in\{4,6\}$, then $A'^{\,d/2}=-I$. 
This will imply the following:
\begin{prop}\label{4or6}
Suppose $A$ is of the form (\ref{a'}), where $A'$ has order $4$ or $6$. 
Then $\semm$ has the \op.
\end{prop}
\begin{proof}
Write $E=\semm$. Let $\p$ be an automorphism and denote by $M$ the matrix representing the induced automorphism $\ba\p$ on $E/\gn{x,t}=\gn{\ba y,\ba z}$.
By the averaging formula, the result will follow once we show that $\R{A^kM}=\infty$ for some $k\in\Z$. Thereto,  
write $\p(t)=z_0t^\e$ for some $\e$ in $\{\pm 1\}$ and  $z_0$ in $\Z^3$. 
A simple calculation shows that $MA'=A'\,^\e M$.
From Section~\ref{realorcomplex}, it immediately follows that $\R M=\infty$ if  $\e=-1$.
So assume $\e=1$. 
Lemma~\ref{com} shows that $M=\pm {A'}^k$ for some $k\in\Z$. 
 Using ${A'}^{d/2}=-I$, we may assume that $M={A'}^k$, 
showing $\R{{A'}^{-k}M}=\R I=\infty$ indeed. 
\end{proof}

So it remains to examine the situation where $A'$ has order~3. 
Tahara \cite[Proposition~3]{t71-1} showed that in this situation, $A$ is conjugate over $\GLL$ to the matrix
\begin{equation}\label{3}
    \begin{pmatrix}
            1&0&\delta\\0&0&-1\\0&1&-1
    \end{pmatrix},
\end{equation}
where $\delta$ is either $0$ or $1$.

\begin{prop}\label{orde3delta1}
Suppose $A$ is of the form (\ref{3}), where $\delta=0$ or $1$. 
Then $\semm$ has \res $6\N\cup\{\infty\}$.
\end{prop}
\begin{proof}
Write $E=\semm$. 
We continue to write 
\[A':=
\begin{pmatrix}
    0&-1\\1&-1
\end{pmatrix}.\]

We first show that $\SE E\subseteq6\N\cup\{\infty\}$. To this end, take an automorphism $\p$ of $E$. 
Write $\ba\p$  the induced automorphism on $E/\gn{x,t^3}$. 
If $\ba\p$ does not induce the identity on $\Z_3$, 
we know from Section~\ref{realorcomplex} that $\R{\ba \p}=\infty$. 
Moreover, $\ZZZ\times \gn{t^3}$ is characteristic in $E$. 
Therefore, if $\R\p$ is finite, we can write
\[\p(x)=x^at^{3b},\; \p(y)=x^cy^dz^et^{3f}, \;\p(z)=x^gy^hz^it^{3j},\; \p(t)=x^ky^lz^mt^n\] for some $a,\lk,n$ in $\Z$, where $n\eq 1$. 
Again, $\p$ is a morphism if and only if it respects the relations \[[y,t]=x^\delta y^{-2}\ii z, \quad[z,t]=\ii z y.\]
Equating $\p(y)\p(z)^{-1}=[\p(z),\p(t)]$ gives
\begin{align*}
x^{c-g}y^{d-h}z^{e-i}t^{3f-3j}&=t^{-3j}z^{-i}y^{-h}x^{-g}\,t^{-n}z^{-m}y^{-l}x^{-k}\,x^gy^hz^it^{3j}\,x^ky^lz^mt^n\\
&=z^{-i}y^{-h}x^{\delta h}y^{-h}z^{-h}y^it^{-n}t^n\\
&=x^{\delta h}y^{-2h+i}z^{-h-i},
\intertext{hence $c=\delta h+g$, $d=-h+i$, $e=-h$ and $f=j$.
Similarly, equating $\p(x)^\delta\p(y)^{-2}\p(z)^{-1}=[\p(y),\p(t)]$ gives}
x^{\delta a -2c-g}y^{-2d-h}z^{-2e-i}t^{3b\delta-6f-3j}&=t^{-3f}z^{-e}y^{-d}x^{-c}\,t^{-n}z^{-m}y^{-l}x^{-k}\,x^cy^dz^et^{3f}\,x^ky^lz^mt^n\\
&=z^{-e}y^{-d}x^{\delta d}y^{-d}z^{-d}y^et^{-n}t^n\\
&=x^{\delta d}y^{-2d+e}z^{-d-e},
\end{align*} 
so we moreover have $\delta a=\delta d+2c+g$ and $\delta b=3f$.

The averaging formula says that $R(\varphi) = \dfrac{1}{3}\R N\sum_{i=0}^2\R{ {A'}^i M}$, with
\[N:=\begin{pmatrix}
    n&b\\3k+\delta(m+l)&a
\end{pmatrix}
\text{ and }
M:=\begin{pmatrix}
    e+i&-e\\e&i
\end{pmatrix}.\]

Since $M$ and $A'$ commute, $M=\pm{A'}^s$ for some $s\in\Z$, hence $\sum_{i=0}^2\R{ {A'}^i M}$ is either $6$ or $\infty$.
Since $n\equiv 1\bmod 3$ and $\delta b=3f$, we moreover have $\R N\in\triples$ as $\delta\in\{0,1\}$. Thus $\R\p\in 6\N\cup\{\infty\}$. 

Conversely, for $\al\in\N$, setting 
\[\p_\al(x)=x^{3\al-1}t^{3(\delta3\al+(1-\delta))},\; \p_\al(y)=x^{\delta\al}y^{-1}t^{3\delta\al}, \;\p_\al(z)=x^{\delta\al}z^{-1}t^{3\delta\al},\; \p_\al(t)=x^{(1-\delta)\al}yt\]
defines an automorphism $\p_\al$ of $E$ inducing the automorphisms
\[N=
\begin{pmatrix}
1&\delta3\al+(1-\delta) \\ \delta+(1-\delta)3\al & 3\al-1
\end{pmatrix}
\text{ and }
M=
\begin{pmatrix}
    -1&0\\0&-1
\end{pmatrix}\in\GL\]
on $\gn{ t^3,x}$ and $\gn{\ba y, \ba z}$, respectively. 
By construction $\R N=3\al$ and $\sum_{i=0}^2\R{ {A'}^i M}=6$, hence $\R\p=6\al$. This completes the proof.
\end{proof}

 \subsubsection{Conclusion}

We summarise the results of this subsection in 
the following table,
where we use the notation $\Lambda(M):=\{\text{eigenvalues of }M\}$. 
 
\[\begin{array}{cc|c}
\hline
\multicolumn{2}{c|}{A}&{\SE \semm}
\\\hline\hline
\rule{0pt}{1em}
1\notin\Lambda(A)
& 
\begin{array}{l}
A=-I\\A\neq -I
\end{array} 
& 
\begin{array}{c}
\even\\\{\infty\}
\end{array} 
\\
\begin{pmatrix}
1&*&*\\0&1&*\\0&0&-1
\end{pmatrix} &  & \{\infty\} 
\\[1em]
\rule{0pt}{3em}
\begin{pmatrix}
1&*&*\\0&-1&n\\0&0&-1
\end{pmatrix} & n\neq 0 & \{\infty\} 
\\[1em]
\rule{0pt}{3em}
\begin{pmatrix}
1&0&\delta\\0&-1&0\\0&0&-1
\end{pmatrix} &
\begin{array}{l}
    \delta=0 \\
     \delta=1
\end{array}&
\begin{array}{c}
    \even \\
     \quadr
\end{array}
\\[1em]
\rule{0pt}{3em}
 \begin{pmatrix}
    1 & \begin{matrix}*&*\,\end{matrix}\\
    \begin{matrix}0\\0\end{matrix} & A'
 \end{pmatrix}&
\begin{array}{l}
    \Lambda(A')\subseteq\real\setminus\{\pm 1\}:\\
    \de {A'}=1 \\
     \de{A'}=-1
\end{array}&
\begin{array}{c}
    \\
    \{\infty\}\text{ or }\{8,\infty\} \\
     \{\infty\}
\end{array}\\[1em]
\rule{0pt}{3em}
 \begin{pmatrix}
    1 & \begin{matrix}*&*\,\end{matrix}\\
    \begin{matrix}0\\0\end{matrix} & A'
 \end{pmatrix}&
\begin{array}{l}
    {A'}^4=I\text{ or }{A'}^6=I \\
     {A'}^3=I\\
\end{array}&
\begin{array}{c}
    \{\infty\} \\
     6\N\cup\{\infty\}\\
\end{array}
\\
\hline
\end{array}
\]

\subsection{\texorpdfstring{The semidirect product $(\ZZ\rtimes_{-I} \Z)\rtimes_\f \Z$}{u}}\label{2by2}

In this subsection, we determine the Reidemeister spectrum of the group $(\ZZ\rtimes_{-I} \Z)\rtimes_\f \Z$. 
Denoting $u$ a generator of the outer $\Z$ and $t$ a generator of the inner $\Z$, the action $\f$ of $u$ on $\ZZ\rtimes_{-I} \Z$ is assumed to be of the form $\f(zt^k)=A(z)(n_0t)^k$, where $n_0\in\ZZ$ and $A$ has either infinite order or $A\neq\pm I$ has order two.

We make the following observation.
\begin{lemma}
The subgroup $\ZZ$ is characteristic in $(\ZZ\rtimes_{-I} \Z)\rtimes_\f \Z$.
\end{lemma}
\begin{proof}
Write $E=(\ZZ\rtimes_{-I} \Z)\rtimes_\f \Z$. 
Let $\p$ be an automorphism of $E$ and take $z$ in $\ZZ$. 
The relation $z^2=ztz\,^{-1}t^{-1}$ shows that $z^2$ belongs to $[E,E]$, which is characteristic. 
Hence $\p(z)^2\in[E,E]$. Moreover, $[E,E]\subseteq\ZZ$ as $E/\ZZ\cong\ZZ$ is abelian.
Hence $\p(z)^2\in\ZZ$, showing $\p(z)\in\ZZ$ as well. 
\end{proof}

\subsubsection{The matrix \ensuremath{A} has order 2}
We start by applying Lemma~\ref{formule} to the situation $A\neq\pm I$ and $A^2=I$.

\begin{prop}\label{Aorder2}
Let $E=(\ZZ\rtimes_{-I} \Z)\rtimes_\f \Z$, where $A=\psi_{|\ZZ}\neq\pm I$ and $A^2=I$. Then $E$ has the \op.
\end{prop}
\begin{proof}
As $A\neq\pm I$ and $A^2=I$, there exists a basis $\{v,w\}$ of $\ZZ$ transforming $A$ into the matrix 
\[
A=\begin{pmatrix}
    1 & r\\
    0 &  -1\\
\end{pmatrix} 
.
\]
Let $\p$ be an automorphism of $E$, and let $\ba \p$ denote the induced automorphism on the quotient $E/\ZZ\cong \ZZ$. 
Write $\ba\p(\ba t)=\ba t^k\ba u^m$,  $\ba\p(\ba u)=\ba t^l\ba u^n$ and $\p(v)=v^aw^c$, $\p(w)=v^bw^d$ so that $\ba \p$ and $\p|_\ZZ$ are represented by the matrices 
\[K=\begin{pmatrix}
    k & l\\
    m & n\\
\end{pmatrix}\text{ and }
M=\begin{pmatrix}
    a & b\\
    c & d\\
\end{pmatrix} \in\GL,\]
respectively.
The action of the quotient $E/\ZZ$ on $\ZZ$ is given by $\al(\ba t^e\ba u^f)=(-I)^eA^f$ for all $e,f\in\Z$.
It is easy to check that $M\al(\ba z)=\al(\ba\p(\ba z))M$ for all $\ba z$ in $E/\ZZ$.
Applying the condition on $\ba t$ gives $-M=(-I)^kA^mM$. 
Hence $k$ is odd and $m$ must be even, which forces $n$ to be odd as $K\in\GL$. 
Applying the condition on $\ba u$ gives 
$MA=(-I)^lA^nM$, so $MA=(-I)^lAM$ since $n$ is odd. 
We distinguish two cases, based on the parity of $l$.

\medskip\noindent\textbf{Even }$\bm{l.}$ 
In this case, $A$ and $M$ commute. Hence $M\in\{\pm I,\pm A\}$ by Lemma~\ref{com}.
At the same time, the \rec $\re{1}{\ba \p}$, $\re{\ba t}{\ba \p}$, $\re{\ba u}{\ba \p}$,  $\re{\ba t\ba u}{\ba \p}$ are all different, as 
\[I-K=\begin{pmatrix}
    1-k & -l\\
    -m & 1-n\\
\end{pmatrix}\in\begin{pmatrix}
    2\Z & 2\Z\\
    2\Z & 2\Z\\
\end{pmatrix}.\]
The addition formula implies that $\R\p$ is at least $\R M+\R{-M}+\R {AM}+\R{-AM}$.
As $M\in\{\pm I,\pm A\}$, one of these terms equals $\R I=\infty$, showing $\R\p=\infty$.

\medskip\noindent\textbf{Odd }$\bm{l.}$ 
Now $MA=-AM$. Equating these matrices explicitly gives
\[\begin{pmatrix}
    a & -b+ar\\
    c & -d+cr\\
\end{pmatrix}=
\begin{pmatrix}
    -a-cr & -b-dr\\
    c & d\\
\end{pmatrix}.
\]
Note that $r=0$ implies $a=d=0$ and $r\neq 0$ implies $a+d=0$. So in both cases $M$ has trace zero. 
The same holds for ${AM}$, for this matrix also satisfies $(AM)A=-A(AM)$.
In addition, the \rec $\re{1}{\ba \p}$ and $\re{\ba u}{\ba \p}$ are different, as 
\[I-K=\begin{pmatrix}
    1-k & -l\\
    -m & 1-n\\
\end{pmatrix}\in\begin{pmatrix}
    \Z & \Z\\
    2\Z & 2\Z\\
\end{pmatrix}.\]
Hence $\R\p$ is at least $\R M+\R {AM}$.
As $M$ and $AM$ both have trace zero, and $\det(A)=-1$, 
either $\di M=1+\det(M)$ or $\di {AM}=1+\det(AM)$ vanishes, thus, either
$\R M$ or $\R {AM}$ is infinite. We conclude that $\R\p=\infty$.
\end{proof}

\subsubsection{The matrix \ensuremath{A} has infinite order}

Similarly, we can determine the \res if $A$ has infinite order. 
We begin by elaborating Lemma~\ref{formule}.

\begin{lemma}\label{formule2}
Let $E=(\ZZ\rtimes_{-I} \Z)\rtimes_\f \Z$, where $A=\psi_{|\ZZ}$ has infinite order. 
Let $\p$ be an automorphism of $E$, write $\p|_\ZZ=M\in\GL$. 
Then $\R\p$, if finite, equals \[\R \p =\R M+\R{-M}+\R {AM}+\R{-AM}.\] 
In addition, either $MA=A^{-1}M$ or $MA=-A^{-1}M$ and the group $\ZZ\rtimes_{-I} \Z$ is characteristic in $E$.
\end{lemma}
\begin{proof}
Let $\ba \p$ denote the induced automorphism on the quotient $E/\ZZ\cong \ZZ$. 
Writing $\ba\p(\ba t)=\ba t^k\ba u^m$,  $\ba\p(\ba u)=\ba t^l\ba u^n$, the automorphism $\ba \p$ is represented by the matrix 
\[K=\begin{pmatrix}
    k & l\\
    m & n\\
\end{pmatrix}\in\GL.
\]
Note that $M\al(\ba z)=\al(\ba\p(\ba z))M$ for all $\ba z$ in $E/\ZZ$, where $\al:\ZZ\to\GL$ is defined by $\ba t^e\ba u^f\mapsto (-I)^eA^f.$
Taking $\ba z=\ba t$ gives $-M=(-I)^kA^mM$.
Hence $m$ must be zero, implying both $k, n\in\{\pm1\}$. Note that this already shows that  $\ZZ\rtimes_{-I} \Z$ is characteristic in $E$.
In case $k=1$ or $n=1$, we have $\R{\ba \p}=\infty$, implying $\R\p=\infty$ as well.
So assume $k=n=-1$. Applying the condition on $\ba u$ gives $MA=(-I)^lA^nM$, hence 
$MA=\pm A^{-1}M$. 
Moreover, the identity 
\[
I-K
=
\begin{pmatrix}
    1-k & -l\\
    -m & 1-n\\
\end{pmatrix}
=
\begin{pmatrix}
    2 & -l\\
    0 & 2\\
\end{pmatrix}
\]
implies that the classes $\re{1}{\ba \p}$, $\re{\ba t}{\ba \p}$, $\re{\ba u}{\ba \p}$ and  $\re{\ba t\ba u}{\ba \p}$ are exactly the \rec of $\ba \p$.
The result now follows from the addition formula.
\end{proof}

If $A$ has infinite order, it has either repeated eigenvalue $\pm 1$ or its eigenvalues are real and different from $\pm 1$. 
In the former case, we have:
\begin{prop}
Let $E=(\ZZ\rtimes_{-I} \Z)\rtimes_\f \Z$, where $A=\psi_{|\ZZ}$ has infinite order. 
If $A$ has repeated eigenvalue $1$ or $-1$, then $E$ has the \op.
\end{prop}
\begin{proof}
Write $\e\in\{\pm1\}$ the repeated eigenvalue of $A$. Changing bases if necessary, we may assume that
$
A=\left(\begin{smallmatrix}
    \e & r\\
    0 &  \e\\
\end{smallmatrix}\right),$
for some $r\neq 0\in\Z.$ 
Clearly $A$ and $-A^{-1}$ have different eigenvalues. 
As a result, they are not conjugate, that is, no $M\in\GL$ satisfies $MA=- A^{-1}M$. 
By Lemma~\ref{formule2}, it suffices to show that $MA=A^{-1}M$ implies $\R M=\infty$ for all $M\in\GL$. 

To this end, take $M$ in $\GL$ satisfying $AMA=M$. Writing
$M=\left(\begin{smallmatrix}
    a & b\\
    c & d\\
\end{smallmatrix}\right),$
this means 
\[
\begin{pmatrix}
    a+\e rc & b+r^2c+\e r(a+d)\\
    c & d+ \e rc\\
\end{pmatrix}
=
\begin{pmatrix}
    a & b\\
    c & d\\
\end{pmatrix}.
\]
Recall $r\neq 0$, hence $c=0$ and subsequently $a+d=0$. 
So either $a$ or $d$ is $1$, implying $\R M=\infty$. This concludes the proof.
\end{proof}

The situation where $A$ has real eigenvalues is more subtle. 
\begin{prop}\label{realeig}
Let $E=(\ZZ\rtimes_{-I} \Z)\rtimes_\f \Z$, where $A=\psi_{|\ZZ}$ has real eigenvalues $\neq \pm1$. 
Then the \res of $E$ is either $\{\infty\}$ or $\{8,\infty\}$. Moreover, if $\sem$ has the \op, then so does $E$.
\end{prop}
\begin{proof}
By Lemma~\ref{formule2}, only $\infty$ and $\R M+\R{-M}+\R {AM}+\R{-AM}$ are candidate \red numbers, where $M\in\GL$ satisfies $MA=\pm A^{-1}M$.

In \cite{gw03-1}, Gon\c calves and Wong  analysed the condition $MA=A^{-1}M$, see Section~\ref{realorcomplex}. They concluded that $MA=A^{-1}M$ can only hold if $\det(A)=1$ and $\T M=0$. 
In particular, $\R M$ is either $2$ or $\infty$, and $\sem$ has the $R_\infty$ property precisely when only $\R M=\infty$ occurs. 
Moreover, if $M$ satisfies the condition $MA=A^{-1}M$, the matrices $-M$, $AM$, $-AM$ satisfy this condition as well. 
Therefore, the sum $\R M+\R{-M}+\R {AM}+\R{-AM}$ is always infinite if $\sem$ has the \op, and either $8$ or $\infty$ otherwise.

We now repeat their argument for the condition $MA=-A^{-1}M$. So, take $M\in\GL$ satisfying $MA=-A^{-1}M$.
Setting $\delta:=\de A$, let $\lm$ and $\delta\ii\lm$ denote the eigenvalues of $A$.
As $\lm\neq\pm1$, these eigenvalues are different. 
In particular $A$ is diagonalizable over $\C$. Take $P$ in $\GC$ such that
\[
PA\ii P
=
\begin{pmatrix}
    \lm & 0 \\
    0 & \delta\ii\lm \\
\end{pmatrix}
\text{ and write }
PM\ii P
=
\begin{pmatrix}
    a & c \\
    b & d \\
\end{pmatrix}.
\]
Denoting $\tilde A:=PA\ii P$ and $\tilde M:=PM\ii P$, we still have the relation $\tilde M\tilde A=-\ii {\tilde A}\tilde M$. 
Equating these matrices explicitly gives
\[
\begin{pmatrix}
    \lm a & \delta\ii\lm c \\
    \lm b & \delta\ii\lm d\\
\end{pmatrix}
=
\begin{pmatrix}
    -\ii\lm a& -\ii\lm c \\
    -\delta \lm b & -\delta\lm d\\
\end{pmatrix}.
\]
As $\lm\neq\pm i$, this implies $a=d=0$. Hence, $\T M=0$. 
As before, this implies $\T {AM}=0$ as well. 
Note that when $\de A=1$, the relation above forces both $b$ and $c$ to be zero, so that $M$ would be the zero matrix. 
Thus $\det(A)=-1$.
In the odd case of Proposition~\ref{Aorder2} we showed that this implies either $\R M=\infty$ or $\R {AM}=\infty.$
In particular, the sum $\R M+\R{-M}+\R {AM}+\R{-AM}$ is always infinite, and the proof is complete.
\end{proof}

Proposition~\ref{realeig} prompts the question: 
\textit{If $\sem$ does not have the \op, does $E$ admit an automorphism having \ren 8?} 
Examples~\ref{wel} and~\ref{niet} below will show that the answer to this question depends on the situation. 

In general, $E$ admits an automorphism having \red number~$8$ if there exists $M\in\SL$, satisfying $MA=A^{-1}M$, that extends to an automorphism $\phi$ of $E$ having finite \red number.
As in fact $G:=\ZZ\rtimes_{-I}\Z$ is characteristic in $E=G\rtimes_\f\Z$ by Lemma~\ref{formule2}, any such $\phi $ is of the form 
\[\phi(gu^k)=\p(g)(g_0u^{-1})^k,\quad g\in G,\,k\in\Z,\]
for some $g_0$ in $G$ and $\p$ in $\Aut G,$ where in turn 
\[\p(zt^k)=M(z)(z_0t^{-1})^k,\quad z\in\ZZ,\,k\in\Z,\]
for some $z_0$ in $\ZZ$.
It is easy to check 
that $\phi$ is an automorphism if and only if $\mu(g_0)\circ\p\circ\f=\ii \f\circ\p$, 
where $\mu(a):G\to G: g\mapsto \con{a}{g}$.
Writing $g_0=m_0t^m$ for some $m_0\in\ZZ$, $m\in\Z$, this implies in particular that  for all $z\in\ZZ$, 
\[\bigl((-1)^mMA\bigr)(z)=\mu(g_0)\bigl(MA(z)\bigr)=\ii AM(z).\]
As $MA=\ii AM$, the above shows that $m$ is even, and that the condition $\mu(g_0)\circ\p\circ\f=\ii \f\circ\p$ is satisfied  on $\ZZ$. 
Hence, $\phi$ is an automorphism if and only if $(\mu(g_0)\circ\p\circ\f)(t)=(\ii \f\circ\p)(t),$ or equivalently, 
$(\mu(m_0)\circ\p\circ\f)(t)=(\ii \f\circ\p)(t)$ as $m$ is even.   
Recalling $\f(t)=n_0t$, this means 
$    (M+\ii A)(n_0)=2m_0+(\ii A-I)(z_0)$, or, equivalently,
 \begin{align}\label{voorw}
    (I+AM)(n_0)&=2A(m_0)+(I-A)(z_0).
\end{align}
So such $\phi$ exists exactly when we can find $m_0, z_0\in\ZZ$ satisfying relation~(\ref{voorw}). 

Clearly, we can find such $m_0$, $z_0$ in the following situations:
\begin{itemize}
    \item If $n_0=0$, that is, $E$ is actually a semidirect product $\ZZ\rtimes\ZZ$.
    \item If $\T A=3$ or $5$, since then $I-A$ is invertible. 
\end{itemize}
A less trivial example where we always can find appropriate $m_0$, $z_0$ is the following: 

\begin{example}\label{wel}
Suppose that $A=\left(\begin{smallmatrix}
    2 & 3\\
    3 & 5\\
\end{smallmatrix}\right),$  and set  $M:=\left(\begin{smallmatrix}
    0 & -1\\
    1 & 0\\
\end{smallmatrix}\right).$ 
It is easy to check that $M\in\SL$ satisfies $AMA=M$, hence $\sem$ does not have the \op.  
Moreover, one easily verifies that
$\im{2A}+\im{I-A}=\ZZ$.
Hence, suitable $m_0, z_0\in\ZZ$ exist for any $n_0\in\ZZ$. 
By the discussion preceding this example, we conclude that $E$ does not have the $R_{\infty}$ property either.
\end{example}

In contrast, the following example shows that the \res of~$E$ may depend on $n_0$. 
\begin{example}\label{niet}
Suppose that $A=\left(\begin{smallmatrix}
    5 & 2\\
    2 & 1\\
\end{smallmatrix}\right).$
The matrix $M$ from Example~\ref{wel} shows that $\sem$ does not have the $R_\infty$ property. 
In general, recall that any $M\in\SL$ satisfying $AMA=M$ can be written as
$
M=
\left(\begin{smallmatrix}
        m&n\\p&-m
\end{smallmatrix}\right),
$ 
where $m$, $n$ and $p$ satisfy system~(\ref{system}). In this particular situation, this system reads: 
\begin{equation*}
\left\{
\begin{aligned}
    &-m^2-pn=1\\
    &4m+2p+2n=0
\end{aligned}
\right.
\end{equation*}

From the bottom equation $-2m=p+n$, it is easy to infer that either both $p$ and $n$ are odd, or both $p$ and $n$ are even. 
Suppose that both $p$ and $n$ are even. Then $m$ is odd by the top equation, so $-m=p/2+n/2$ moreover implies that $p$ or $n$ is a quadruple. 
Consequently, $pn\equiv 0\bmod 8$, thus $-1\equiv m^2 \bmod 8$ is a square modulo~$8$. 
This is false, hence $p$ and $n$ must be odd, implying $m$ to be even. We conclude that
all entries of $I+M$ are odd.
Of course, 
this
holds equally for the matrix $I+AM$, as $AM$ satisfies $A(AM)A=AM$ as well. 
As all entries of $I-A$ are even, 
equation~(\ref{voorw}) allows a solution $(m_0,z_0)$ precisely when the sum of the components of $n_0$ is even. 
Equivalently, $E$ has the $R_\infty$ property precisely when the sum of the components of $n_0$ is odd. 
\end{example}

\subsubsection{Conclusion}

We summarise our findings in the following table:
\[
\begin{array}{c|c}
\hline
\text{eigenvalues of }A:=\f|_\ZZ &\SE {(\ZZ\rtimes_{-I}\Z)\rtimes_\f\Z}\\ \hline\hline
\rule{0pt}{1em}
 \begin{array}{c}
1,1\text{ or }-1,-1\\
\lm, -\ii \lm,\; \lm\in\real\\
\lm, \ii \lm,\; \lm\in\real\setminus\{\pm1\}
\end{array}&
\begin{array}{c} 
\{\infty\}\\ 
\{\infty\}\\
\{\infty\}\text{ or }\{8,\infty\}
\end{array}\\
\hline
\end{array}
\] 

\subsection{\texorpdfstring{The semidirect product $H_n\rtimes_\f \Z$}{u}}\label{delaatstereferentie}

Let us first fix some notations we will use throughout this subsection.
As before, we will write elements of ${H_n}$ as expressions of the form $x^ay^bz^c$, 
where $z$ generates the center of $H_n$ and $[y,x]=z^n$.  
Denote the projection map $H_n\to H_n/\zn$ by $h\mapsto \ba h$. 
The quotient $H_n/\zn$ is generated by $\ba x$ and $\ba y$, and moreover $\gn {\ba x, \ba y}\cong\ZZ$. 
Let $A\in\GL$ represent the induced automorphism $\ba\f$ relative to the basis $\{\ba x, \ba y\}$.

Again, we have:
\begin{lemma}\label{charhn}
If $A$ does not have $1$ as eigenvalue, the subgroup $H_n$ of $\sehm$ is characteristic. 
\end{lemma}
\begin{proof}
Write $E=\sehm$. 
Since $A$ does not have $1$ as eigenvalue, 
$\im{I-A}$ has finite index in $H_n/\zn$. Moreover, 
one easily checks that
$\im{I-A}= \ov\ab$, . 
so $\ov\ab$ has finite index in $H_n/\zn$ as well.
Hence $\ab$ has finite index in $H_n$ and the result follows.
\end{proof}

In fact, if $A$ does have eigenvalue $1$, we already studied the \res of $\sehm$. Indeed, if $A$ has repeated eigenvalue $1$, the group $\sehm$ is nilpotent; if $A$ has eigenvalues $1$ and $-1$, the group $\sehm$ is an extension of $\ZZ$ by $\ZZ$.
\begin{lemma}
If $A$ has eigenvalues $1$ and $-1$, the group $\sehm$ is an extension of $\ZZ$ by $\ZZ$.
\end{lemma}
\begin{proof}
Write $E=\sehm$. 
The matrix $I-A$ has eigenvalues $0$ and $2$.
Consequently, $\im{I-A}\cong \Z$, say $\im{I-A}=\gn {\ba h}$ for some $h\in H_n$.
Take  generators $\ba v$, $\ba w$ of $H_n/\zn$ such that ${\ba v}^k=\ba h$ for some $k\in\Z$, 
or equivalently, $\im{I-A}\subseteq \gn{\ba v}$.
Since $\ov{[E,E]}= \im{I-A}$, 
also  $\ov{[E,E]}\subseteq \gn{\ba v}$, that is, $[E,E]\subseteq \gn{v,z}$. 
It easily follows that $\gn{v,z}\nor E$ and that $E/\gn{v,z} 
\cong\ZZ$.
Clearly, also  $\gn{v,z}\cong\ZZ$, and the proof is finished.
\end{proof}
So it remains to study $\sehm$ when $A$ does not have $1$ as eigenvalue.

\begin{prop}
If $A\neq -I$ does not have $1$ as eigenvalue, the group $\sehm$ has the \op.
\end{prop}
\begin{proof}
The center $\zn$ is characteristic in $\sehm $ and
$(\sehm)/\zn\cong\sem$. 
Therefore, if $\sem$ has the \op, $\sehm$ has the $R_\infty$ property as well.
It thus remains to prove the proposition when $A$ has real eigenvalues and $\de A=1$. 
To this end, take an automorphism $\phi$ of $\sehm$. 
If $\phi$ induces the identity on $(\sehm)/H_n\cong\Z$, clearly $\R\phi$ is infinite. 
So assume $\phi(t)=h_0\ii t$ for some $h_0\in H_n$.
Let $\p$ denote the restriction of $\phi$ to $H_n$.
Further, let $\ba\p$ denote the induced automorphism on the quotient $H_n/\zn$,
say $\ba\p$ is represented by $M\in\GL$. 
The addition formula implies that $\R\p=\R{\ba \p}\R{\p|_{\zn}}$.

From the condition $\phi(th)=\phi(t)\phi(h)$ for all $h\in H_n$, it easily follows that
$MA=\ii A M$. 
Hence
$\R{\ba\p}=\infty$ when $\de M=-1$.

At the same time, $\p|_{\zn}$ is multiplication by $\de M$, hence $\R{\p|_{\zn}}=\infty$ when $\de M=1$. 
We conclude that $\R\p$ is always infinite, 
so $\R\phi=\infty$ as well.
\end{proof}

In contrast, when $A=-I$, more is possible:

\begin{prop}
Let $E$ denote the group $\sehm$ where $\f(x)=\ii xz^k$ and $\f(y)=\ii y z^l$ for some $k,l\in\Z$.
\begin{enumerate}
    \item If $n$ is odd, $E$ has \res $\quadr$.
    \item If $n$ is even, $E$ has \res $\quadr$ if both $k$  and $l$ are  even, and \res $8\N\cup\{\infty\}$ if $k$ or $l$ is odd.
\end{enumerate}
\end{prop}
\begin{proof}
We first show that in either case, the \res of $E$ is contained in $4\N\cup\{\infty\}$. 
To this end, let $\phi$ be an automorphism of $E$. 
We may assume that $\phi(t)=h_0\ii t$ for some $h_0\in H_n$.
As before, write $\p$ the restriction to $H_n$ 
and let $\ba \p$ denote the induced automorphism on the quotient $H_n/\zn$.
Further, let $M\in\GL$ denote the matrix representing $\ba\p$ relative to the basis $\{\ba x, \ba y\}$.
If $\de M=1$, then $\p$ induces the identity on $Z(H_n)$, hence $\R\phi=\infty$. 
Furthermore, if $\de M=-1$, the addition formula says that $\R\phi=4\,|\T M|$ if $\T M$ is nonzero, 
and $\R\phi=\infty$ otherwise.
So $\R\p\in\quadr.$

Next, we determine which $M$ arise as induced automorphisms on $H_n/\zn$, or more generally, which $\p$ can be extended to an automorphism on $E$. 
As elements in $\SL $ always contribute $\infty$ to $\SE E$, we only consider $M$ having $\de M=-1$.
So, let $\p\in\Aut {H_n}$ and write $\p(x)=x^ay^cz^m$, $\p(y)=x^by^dz^p$ for some $m,p\in\Z$ and 
\begin{equation}\label{not}
    M:=
    \begin{pmatrix} a & b\\ c & d\end{pmatrix}\in\GL, 
    \text{ where }\de M=-1.
\end{equation}
Then $\p$ extends to an automorphism $\phi:E\to E:ht^r\mapsto \p(h)(h_0\ii t)^r$, $h_0\in H_n$,  if and only if all $h\in H_n$ satisfy 
$\phi(\f(h)t)=\phi(t)\phi(h)$, that is,
\[(\p\ci\f)(h)\,h_0=h_0(\ii\f\ci\p)(h).\]
Writing $h_0=x^e y^fz^g$ (where $g$ will not really matter as $z^g$ is central), a simple calculation shows this to be equivalent to
\begin{equation}\label{con}
    2\begin{pmatrix}m\\p\end{pmatrix}+(I+M^T)\begin{pmatrix}k\\l\end{pmatrix}
    =
    n\left(
    \begin{pmatrix}ac\\bd\end{pmatrix}
    +
    \begin{pmatrix}
         -c & a\\
         -d & b
    \end{pmatrix}
    \begin{pmatrix}e\\f\end{pmatrix}
    \right).
\end{equation}

If $n$ is odd, any $M\in\GL$ allows $m,p, e, f\in\Z$ such that condition~(\ref{con}) above is satisfied. 
Indeed, since $n$ is odd, $2$ is invertible modulo $n$, 
so there exist $m,p\in\Z$ such that 
\[
2\begin{pmatrix}m\\p\end{pmatrix}+(I+M^T)\begin{pmatrix}k\\l\end{pmatrix}\in\begin{pmatrix}n\Z \\ n\Z\end{pmatrix},\text{ say }\,
2\begin{pmatrix}m\\p\end{pmatrix}+(I+M^T)\begin{pmatrix}k\\l\end{pmatrix}=\begin{pmatrix}n\ti m \\ n\ti p\end{pmatrix}
\]
for some $\ti m, \ti p\in\Z$.
Subsequently, setting
\[
\begin{pmatrix}e\\ f\end{pmatrix}
:= 
\ii{\begin{pmatrix}
         -c & a\\
         -d & b
    \end{pmatrix}} 
\begin{pmatrix}\ti m-ac\\ \ti p -bd\end{pmatrix}
\] results in a solution $(m,p,e,f)$ to (\ref{con}).
Therefore, there exists for every matrix 
$M_r:=\left(\begin{smallmatrix} r &1\\1&0\end{smallmatrix}\right),$ $r\in\N,$ 
an automorphism $\phi_r$ of $E$ inducing $M_r$ on $H_n/\zn$. 
By construction $\R{\phi_r}=4r$, thus $E$ has \res $4\N\cup\{\infty\}$. 

If $n$ is even, we can similarly find suitable $m, p, e, f\in\Z$ provided 
\begin{equation}\label{even}
    (I+M^T)\begin{pmatrix}k\\l\end{pmatrix}\in\begin{pmatrix}2\Z \\ 2\Z\end{pmatrix},
\end{equation} 
to guarantee the existence of $m$ and $p$. 
Moreover, condition~(\ref{con}) readily implies that  if $n$ is even, (\ref{even}) is in fact equivalent to the existence of such $m$, $p$, $e$ and $f$. 
Hence $M$ lifts to an automorphism of $E$ if and only if it satisfies condition~(\ref{even}) above.
Therefore, the \res of $E$ equals the set 
\[\SE E=\{4|\T M|\mid M\text{ satisfies (\ref{even})}, \de{M}=-1, \T M\neq 0\}\cup\{\infty\}.\]
We distinguish four cases, based on the parity of $k$ and $l$.

\medskip\noindent\textbf{Case 1.}
If both $k$ and $l$ are even, condition~(\ref{even}) is trivially satisfied. 
As before, the matrices $M_r:=\left(\begin{smallmatrix} r &1\\1&0\end{smallmatrix}\right),$ $r\in\N,$ show that $E$ has \res $4\N\cup\{\infty\}$.

\medskip\noindent\textbf{Case 2.} Next, suppose that $k$ is even, but $l$ is odd. 
Let $M\in\GL$ satisfy condition~(\ref{even}), write $M$ as in (\ref{not}).
Condition~(\ref{even}) implies that both $c$ and $1+d$ are even, that is, $c$ is even and $d$ is odd. 
As $M$ is invertible, $a$ must be odd as well. 
Hence $\T M=a+d$ is even, showing $\SE E\subseteq 8\N\cup\{\infty\}$.

Conversely, for $r\in\N$ even, the matrix 
$M_r:=\left(\begin{smallmatrix}
    r+1 & r/2 \\ -2 & -1
\end{smallmatrix}\right)$
satisfies condition~(\ref{even}), has $\de {M_r}=-1$ and moreover $\T{M_r}=r$. 
Thus in this case, indeed $\SE E= 8\N\cup\{\infty\}$.

\medskip\noindent\textbf{Case 3.} Note that $x$ and $y$ play quite a symmetric role: 
we only use that $H_n/\zn=\gn{\ba x, \ba y}$ and that $[y,x]=z^n$, where $\zn=\gn z$. 
Consequently, if $k$ is odd and $l$ is even, 
we can reduce to case~2 by swapping $x$ and $y$ and replacing $z$ by $\ii z$.

\medskip\noindent\textbf{Case 4.} Similarly, if both $k$ and $l$ are odd, replacing $x$ by $yx$ transforms this case into case~2, 
as $\f(yx)=(yx)^{-1}z^{k+l+n}$ and $n$ is even, and $[y,yx]=[y,x]$. 
\end{proof}

\subsubsection{Conclusion}

We summarise our findings in 
the following table:
\medskip
\[\begin{array}{ll|c}
\hline
\multicolumn{2}{c|}{\f}&\SE {H_n\rtimes_\f\Z} \\
\hline\hline
\rule{0pt}{4ex}
A:=\f|_{\ZZ}\neq -I\text{ and }1\notin\Lambda(A)& &\{\infty\}\\[1em]
\f(x)=\ii xz^k
, 
\f(y)=\ii yz^l
&
\begin{array}{c}(k\text{ and }l\text{ even) or }(n\text{ odd})\\
(k\text{ or }l\text{ odd) and }(n\text{ even})
\end{array}
&
\begin{array}{c}4\N\cup\{\infty\}\\
8\N\cup\{\infty\}\end{array}\\[4ex]
\hline
\end{array}
\]

\printbibliography

\end{document}